\documentclass[amsthm]{amsart}

\usepackage{amscd}
\usepackage{amsmath,amsthm}
\usepackage{amssymb}
\usepackage{units}
\usepackage{listings}
\usepackage{enumitem}
\usepackage{emptypage}
\usepackage{tikz}
\usetikzlibrary{shapes}
\usetikzlibrary{arrows}
\usepackage[all]{xy}
\usepackage{color}
\usepackage{algorithm}
\usepackage[noend]{algcompatible}
\usepackage{fancyhdr}
\usepackage{hyperref}
\hypersetup{pdftex,colorlinks=true,allcolors=black}
\usepackage{hypcap}
\usepackage{stmaryrd}
\usepackage{etoolbox}
\usepackage{listings}
\usepackage{calc}

\def\nn{\mathbb N}              
\def\qq{{\mathbb Q}}
\def\zz{{\mathbb Z}}

\def\opn#1#2{\def#1{\operatorname{#2}}} 
\opn\chara{char} \opn\length{\ell} \opn\pd{pd} \opn\rk{rk}
\opn\projdim{proj\,dim} \opn\injdim{inj\,dim} \opn\rank{rank}
\opn\depth{depth} \opn\sdepth{sdepth} \opn\hdepth{hdepth}
\opn\grade{grade} \opn\height{height} \opn\embdim{emb\,dim}
\opn\codim{codim}  
\opn\Tr{Tr} \opn\bigrank{big\,rank}
\opn\superheight{superheight}
\opn\lcm{lcm}
\opn\trdeg{tr\,deg}
\opn\reg{reg} \opn\lreg{lreg} \opn\ini{in} \opn\lpd{lpd}
\opn\size{size}

\opn\div{div} \opn\Div{Div} \opn\cl{cl} \opn\Cl{Cl}
\opn\Spec{Spec} \opn\Max {Max} \opn\Supp{Supp} \opn\supp{supp} \opn\Sing{Sing}
\opn\Ass{Ass} \opn\Min{Min} \opn\rang{rang}
\opn\hait{ht}
\opn\Ann{Ann} \opn\Rad{Rad} \opn\Soc{Soc}
\opn\Im{Im} \opn\Ker{Ker} \opn\Coker{Coker} \opn\Am{Am}
\opn\Hom{Hom} \opn\Tor{Tor} \opn\Ext{Ext} \opn\End{End}
\opn\Aut{Aut}
\opn\id{id}  
\opn\deg{deg}
\opn\mod{mod}
\opn\ord{ord}
\opn\Syz{Syz}
\opn\Hilb{Hilb}
\opn\aff{aff} \opn\con{conv} \opn\relint{relint} \opn\st{st}
\opn\lk{lk} \opn\cn{cn} \opn\core{core} \opn\vol{vol}
\opn\link{link} \opn\star{star}
\opn\gr{gr}

\def\pot#1#2{#1[\kern-0.28ex[#2]\kern-0.28ex]}

\opn\Mon{Mon}
\opn\lm{LM}
\opn\lexp{LE}
\opn\lt{LT}
\opn\sig{sig}
\opn\lc{LC}
\opn\z{\mathbb Z}
\opn\q{\mathbb Q}
\opn\tail{tail}
\opn\quot{Quot}
\opn\L{L}
\opn\nf{NF}
\opn\red{red}
\opn\ecart{ecart}
\opn\sp{s-poly}
\opn\smiss{smiss}
\opn\gmiss{gmiss}
\opn\syz{Syz}
\opn\hc{HC}
\opn\h{H}
\opn\hp{HP}
\opn\p{P}
\opn\gp{gcd-poly}
\opn\ep{ext-poly}
\opn\ssyz{\sp_{syz}}
\opn\gsyz{\gp_{syz}}
\def\sing{\textsc{Singular }}
\opn\heit{ht}

\opn\dirlim{\underrightarrow{\lim}}
\opn\inivlim{\underleftarrow{\lim}}

\let\Dirsum=\bigoplus

\newcommand{\gen}[1]{\left\langle #1 \right\rangle}

\newcommand{\bpa}[1]{\big( #1 \big)}

\newcommand{\set}[2]{\left\{ #1\ \mid\ #2\right\}}
\newcommand{\aco}[1]{\left\{ #1 \right\}}
\newcommand{\itce}{\item[$\circ$]}

\newcommand{\fracs}[2]{\displaystyle\frac{#1}{#2}}
\newcommand{\randex}[2]{$\arraycolsep=0pt\begin{array}{lr}\textnormal{{#1}:}&\textnormal{#2}\end{array}$}
\newcommand{\Dsum}[2]{\displaystyle\Dirsum_{#1}^{#2}}
\newcommand{\Sum}[2]{\displaystyle\sum_{#1}^{#2}}

\newcommand{\tinm}[1]{\textnormal{\tiny{#1}}}
\usepackage{color}

\newcommand{\hide}[1]{}

\lstnewenvironment{linux}{\lstset{language=bash}\lstset{
 morekeywords={git,gdb}
}}{}
\lstnewenvironment{singular}{\lstset{language=C++}\lstset{
 morekeywords={poly,ring,ideal,reduce,proc,number,jet,sumcoef,deg,inverse,system,dim,LIB,desingularization,
 factorial,invp,string,var}}}{}
\lstnewenvironment{cpp}{\lstset{language=C++}\lstset{
 morekeywords={poly,ring,ideal}
}}{}
\let\epsilon\varepsilon
\let\phi=\varphi
\let\kappa=\varkappa

\def\qed{\ifhmode\textqed\fi
      \ifmmode\ifinner\quad\qedsymbol\else\dispqed\fi\fi}
\def\textqed{\unskip\nobreak\penalty50
       \hskip2em\hbox{}\nobreak\hfil\qedsymbol
       \parfillskip=0pt \finalhyphendemerits=0}
\def\dispqed{\rlap{\qquad\qedsymbol}}

\opn\dis{dis}
\def\pnt{{\raise0.5mm\hbox{\large\bf.}}}

\opn\Lex{Lex}

\def\ret{\textbf{return}}
\newtheorem{alg}{Algorithm}

\definecolor{codegreen}{rgb}{0,0.6,0}
\definecolor{codegray}{rgb}{0.5,0.5,0.5}
\definecolor{codepurple}{rgb}{0.58,0,0.82}
\definecolor{backcolour}{rgb}{0.99,0.99,0.99}
\lstdefinestyle{adistyle}{
    backgroundcolor=\color{backcolour},   
    commentstyle=\color{codegreen},
    keywordstyle=\color{blue},
    numberstyle=\tiny\color{codegray},
    stringstyle=\color{codepurple},
    basicstyle=\footnotesize\ttfamily,
    breakatwhitespace=false,         
    breaklines=true,                 
    captionpos=b,                    
    keepspaces=true,                 
    numbers=none,                    
    showspaces=false,                
    showstringspaces=false,
    showtabs=false,                  
    tabsize=2,
    xleftmargin=0.5cm
}
\newtheorem{defn}{Definition}
\newtheorem{ex}{Example}

\newtheorem{rem}{Remark}

\begin{document}

\title{New Strategies for Standard Bases over Z}


\author{Christian Eder}
\address{Department of Mathematics, University of Kaiserslautern,\\ Erwin-Schr\"odinger-Str., 67663 Kaiserslautern, Germany}
\email{ederc@mathematik.uni-kl.de}

\author{Gerhard Pfister}
\address{Department of Mathematics, University of Kaiserslautern,\\ Erwin-Schr\"odinger-Str., 67663 Kaiserslautern, Germany}
\email{pfister@mathematik.uni-kl.de}

\author{Adrian Popescu}
\address{Department of Mathematics, University of Kaiserslautern,\\ Erwin-Schr\"odinger-Str., 67663 Kaiserslautern, Germany}
\email{popescu@mathematik.uni-kl.de}

\begin{abstract}
Experiences with the implementation of strong Gr\"obner bases respectively standard bases for polynomial rings over principal ideal rings are explained: different strategies for creating the pair set, methods to avoid coefficient growth and a normal form algorithm for non-global orderings.
\end{abstract}


\maketitle

\section{Introduction}

The aim of this paper is to explain our experience in implementing the strong Gr\"obner bases (resp. standard bases) algorithm for polynomial rings over principal ideal rings. We use the following notation (for more details see \cite{singular-book}).

\begin{defn}
Let $\textnormal{Mon}(x_1, \ldots, x_n)$ be the set of monomials in $n \in \nn$ variables $x_1, \ldots, x_n$. We denote by $x$ the set of all variables and for $\alpha = (\alpha_1,\ldots,\alpha_n) \in \nn^n$ we set $x^\alpha := x_1^{\alpha_1}\cdot \ldots \cdot x_n^{\alpha_n}$. 

A \textbf{monomial ordering} $<$ is a total ordering on $\textnormal{Mon}(x)$ satisfying the following property $$x^\alpha < x^\beta \implies x^\gamma \cdot x^\alpha < x^\gamma\cdot x^\beta \quad \textnormal{for all } \alpha,\beta,\gamma \in \nn^n.$$

A monomial ordering is called \textbf{global} if $1 = x^0 \leq x^\alpha$ for all $\alpha \in \nn^n$.

A monomial ordering is called \textbf{local} if $1 = x^0 \geq x^\alpha$ for all $\alpha \in \nn^n$.

A monomial ordering is called \textbf{mixed} if there is $\alpha,\beta \in \nn^n$ such that $x^\alpha < 1 < x^\beta$.
\end{defn}

\begin{defn}\label{def: leading stuff} Let $<$ be a monomial ordering on Mon$(x)$, $R$ be a ring and $f \in R[x]:=R[x_1, \ldots, x_n]$ a polynomial.  We can write $f$ in a unique way $$f = \Sum{i=1}{L}c_i \cdot x^{\alpha^{(i)}},$$ $0 \neq c_i \in R$ and $\alpha^{(i)} \in \nn^n$ for $1\leq i \leq L$ such that $x^{\alpha^{(1)}} < \ldots < x^{\alpha^{(L)}}$. We define
\begin{itemize}
\itce the \textbf{length} of $f$ denoted by $\length(f) = L$,
\itce the \textbf{leading term} denoted by $\lt(f) = c_L\cdot x^{\alpha^{(L)}}$,
\itce the \textbf{leading coefficient} denoted by $\lc(f) = c_L$,
\itce the \textbf{leading monomial} denoted by $\lm(f) = x^{\alpha^{(L)}}$,
\itce the \textbf{tail} of $f$ denoted by $\tail(f) = f - \lt(f)$,
\itce the \textbf{degree} of $f$ denoted by $\deg(f) = \displaystyle\max_{1\leq i \leq L}\left\{\alpha^{(i)}_1 + \ldots + \alpha^{(i)}_n\right\}$ and
\itce the \textbf{ecart} of $f$ denoted by $\ecart(f) = \deg(f) - \deg\big(\lt(f)\big)$.
\end{itemize}
\end{defn}

\begin{defn} Let $I \subseteq R[x]$ be an ideal. 

The \textbf{leading ideal} of $I$ is L$_<(I) = \gen{\set{\lt(f)}{f\in I}}$. 

The \textbf{leading monomial ideal} of $I$ is LM$_<(I) = \gen{\set{\lm(f)}{f\in I}}$. 

The leading ideal for a set of polynomials $S:= \aco{f_1,\ldots, f_m}$ is defined to be the leading ideal of the ideal generated by $S$.

From now on, if the monomial ordering is clear from the context, we will only write L$(I)$ (respectively LM$(I)$) instead of L$_<(I)$ (respectively LM$_<(I)$).
\end{defn}



\begin{defn}[Standard Bases]\label{def: standard bases}
A \textbf{standard basis} of an ideal $I \subset R[x]$ with respect to a fixed monomial ordering $<$ is a finite set $S \subset I$ such that L$(I) =$ L$(S)$. 

$S$ is a \textbf{strong standard basis} if additionally it satisfies the following property: $$\forall f \in I , \exists\ g \in S \textnormal{ such that } \lt(g)\ |\ \lt(f).$$
\end{defn}

Next we present a simple Buchberger algorithm for computing standard bases in polynomial rings over principal ideal rings. For this we introduce some more basic notions.

In the following, let $R$ be a ring. Usually we will take $R$ to be $\zz$ or $\zz_r := \nicefrac{\zz}{r\zz}$ where $r\in \zz$ is not a prime number. Let $f,g\in R[x]$ be two polynomials such that $$f = c_f\cdot m_f + \ldots$$ $$g = c_g\cdot m_g + \ldots$$ with $\lt(f) = c_f \cdot m_f$ and $\lt(g) = c_g\cdot m_g$.

\begin{defn}\label{def: spair,gcdpair,extpoly} We define several important polynomials used in the strong standard bases theory.
\begin{itemize}
\itce The \textbf{$s-$polynomial} (or \textbf{s$-$pair}) of $f$ and $g$, denoted by $$\sp(f,g) = \fracs{\lcm(c_f,c_g)\lcm\left(m_f,m_g\right)}{c_f m_f}\cdot f - \fracs{\lcm(c_f,c_g)\lcm\left(m_f,m_g\right)}{c_g m_g}\cdot g.$$ Note that the leading terms will cancel out. 
\itce Since $R$ is a principal ideal ring, we have that $\gen{c_f,c_g} = \gen{c}$. Let $c = d_f \cdot c_f + d_g \cdot c_g$.
The \textbf{strong polynomial} (or the \textbf{gcd polynomial}, \textbf{gcd$-$pair}, \textbf{strong pair}) of $f$ and $g$ denoted by $$\gp(f,g) = \fracs{d_f \lcm(m_f,m_g)}{m_f} \cdot f + \fracs{d_g \lcm(m_f,m_g)}{m_g} \cdot g.$$ Notice that $\lt\bpa{\gp(f,g)} = c \cdot \lcm(m_f,m_g)$.
\itce In case of rings with zero divisors, the \textbf{extended s$-$polynomial} of $f$ denoted by $$\ep(f) = \Ann(c_f) \cdot f = \Ann(c_f) \cdot \tail(f).$$
Note that if $\lc(f)$ is not a zero divisor then $\ep(f) = 0$.
\end{itemize}
\end{defn}

\begin{rem}Note that in the definition of gcd$-$pairs $d_f$ and $d_g$ are not uniquely determined. 
%
%
Let two polynomials $f,g \in R[x]$, $$f = c_f \cdot m_f + \ldots,$$ $$g = c_g \cdot m_g + \ldots,$$ such that $\gen{c} = \gen{c_f,c_g}$. Let $d_f \neq d'_f,d_g \neq d'_g$ such that $$d_f \cdot c_f + d_g \cdot c_g = c = d'_f \cdot c_F + d'_g \cdot c_g.$$
Denote by $P$ and $Q$ the gcd$-$polynomial constructed with $d_f$,$d_g$, respectively $d'_f$ and $d'_g$. Then $P = Q + \alpha \cdot  \sp(f,g)$ for an $\alpha \in R$.

\end{rem}

%


Note that if $\gen{c_f,c_g} = \gen{c_f}$, then we can take $d_f = 1$ and $d_g = 0$. Therefore \linebreak $\gp{(f,g)} = m \cdot f$, for the monomial $m = \fracs{\lcm(m_f,m_g)}{m_f}$. In other words if one of the leading coefficients divides the other, then the $\gp$ is equal to a multiple of the initial polynomial and, as we will later see, it brings no new information. We call such a gcd$-$polynomial \textbf{redundant}.

\begin{alg}[Buchberger's algorithm]\label{alg: std blueprint} Let $R$ be a principal ideal ring, $<$ a monomial ordering on $R[x]$.

\begin{algorithm}[H]
\caption{std(ideal $I$)}
\begin{algorithmic}[1]
\REQUIRE $I = \gen{f_1, \ldots, f_r} \subseteq R[x]$ an ideal
\ENSURE $S = \aco{g_1, \ldots, g_t}$ is a strong standard basis of $I$
\STATE $S = \{f_1, \ldots, f_r\}$;
\STATE $L = \set{\sp(f_i, f_j)}{i < j} \cup \set{\gp(f_i, f_j)}{i < j} \cup \aco{\ep(f_i)}$;
\WHILE{$L \neq \emptyset$}
\STATE choose $h \in L$ and \footnotemark reduce   it with $S$;
\IF{$h \neq 0$}
\STATE $S = S \cup \aco{h}$;
\STATE $L = L \cup \{ \sp(g,h) \ | \ g \in S \} \cup \set{\gp(g,h)}{g \in S} \cup \aco{\ep(h)}$;
\ENDIF
\ENDWHILE
\ret\ $S$
\end{algorithmic}
\end{algorithm}
\footnotetext{Normal form of $h$ with respect to $S$}

\end{alg}

In the second section we compare two implementations of Buchberger's algorithm. They are distinguished in the modifications of step 7 in Algorithm \ref{alg: std blueprint}.

In the third section we explain a method to avoid coefficient growth.

The fourth section describes the normal form procedure for local orderings which is implemented in the standard basis algorithm. The normal form procedure used for polynomial rings over fields does not terminate for local orderings over rings.

\section{ALL vs. JUST}\label{section: all vs just}

In this section let $R = \zz$. 
D. Lichtblau proved in \cite[Theorem 2]{LB} that in Algorithm \ref{alg: std blueprint} it is in fact enough to take for a pair $(g,h)$ either the s-poly or the gcd-poly. Note line 7 of the pseudo-code where we add $\sp(g,h)$ and $\gp(g,h)$ to the pair list $L$. Instead of this, if $\gp(g,h)$ is not redundant then we only add $\gp(g,h)$ and if $\gp(g,h)$ is redundant then we only add $\sp(g,h)$. 

We name the usual strategy in which we consider all pairs by ALL and the one in which we consider just one pair simply by JUST. 

At a first glance fewer pairs could be interpreted as a faster computation and less used memory. This turns out to be wrong, since the few pairs from JUST will generate many pairs later in the algorithm. We have compared the timings of the two different strategies for random examples. 

In most cases, ALL strategy was faster than JUST. We found examples for which ALL was 38\,000 times faster than JUST and examples where JUST was 2\,300 times faster than ALL.
In the tables below we present the different timings we obtained for our examples. All of the timings are represented in milliseconds. The input ideals can be found in the Appendix 
.

\begin{table}[H]
\begin{center}

\begin{tabular}{l||r|r|r|}

 & \multicolumn{1}{c|}{ALL time} &  \multicolumn{1}{c|}{JUST time} &  \multicolumn{1}{c|}{Factor}  \\ \hline
Example \randex{A}{1} \hide{2:106} &  6\,630 & 251\,481\,340 & 38\,000 \\ \hline

Example \randex{A}{2} \hide{1:250} &  40 & 63\,740 & 1\,600 \\ \hline

Example \randex{A}{3} \hide{2:109} & 51\,570 & 6\,134\,060 & 119  \\ \hline

Example \randex{A}{4} \hide{1:238} & 4\,400 & 522\,340 & 118 \\ \hline

Example \randex{A}{5} \hide{1:227} & 250 & 19\,080 & 76 \\ \hline 

Example \randex{A}{6} \hide{1:265} & 4\,110 & 271\,460 & 66 \\ \hline 

Example \randex{A}{7} \hide{2:125} & 56\,720 & 2\,676\,080 & 47 \\ \hline

Example \randex{A}{8} \hide{2:105} & 190 & 8\,340 & 44 \\ \hline 

Example \randex{A}{9} \hide{2:141} & 1\,131\,470 & 39\,694\,950 & 35 \\ \hline 

Example \randex{A}{10} \hide{0:66} & 110 & 3\,020 & 27 \\ \hline

\end{tabular}
\caption{Best examples for ALL strategy}
\label{table: ALL}
\end{center}
\end{table}

\begin{table}[H]
\begin{center}

\begin{tabular}{l||r|r|r|}
 &  \multicolumn{1}{c|}{ALL time} &  \multicolumn{1}{c|}{JUST time} &  \multicolumn{1}{c|}{Factor}  \\ \hline
 
Example \randex{B}{1} \hide{3:1} & 1\,696\,700 & 730 & 2\,324 \\ \hline

Example \randex{B}{2} \hide{0:79} &  252\,950 & 450 & 562 \\ \hline

Example \randex{B}{3} \hide{3:8} & 4\,090 & 30 & 136 \\ \hline

Example \randex{B}{4} \hide{3:6} & 35\,160 & 500 & 70 \\ \hline

Example \randex{B}{5} \hide{3:12} & 3\,690 & 110 & 33 \\ \hline

\end{tabular}
\caption{Best examples for JUST strategy}
\label{table: JUST}
\end{center}
\end{table}

All these examples where randomly generated. We considered the ideals over $\zz[x,y,z]$ with the degree reverse lexicographical ordering.  
%
%
%
%

We searched  for such examples over finite rings, like $\zz_{2^{100}}$ or $\zz_{10^{200}}$ but we were unable to find interesting examples over these finite rings. The reason why we found examples over $\zz$ is that here there are no bounds for coefficients. In order to give an idea of how huge the coefficients can be in the standard bases computation over $\zz$, we computed Example B:1 over  $\zz_{10^{200}}$ with the ALL strategy and got the result after 8 seconds; over $\zz_{10^{1000}}$ we obtained the result after 648 seconds (10 minutes) and over $\zz$ we got the result after almost 30 hours. One can imagine how big the coefficients are when computing over $\zz[x]$. This is the main reason why in some cases, the computation of a standard basis over $\zz$ will be slow.

\section{Huge coefficients over the integers}\label{sec: huge coeffs}

One of the most common problems when computing over $\zz[x]$ is that the coefficients will rapidly increase. If during the algorithm we have two polynomials with big coefficients, when computing their s$-$polynomial (or strong polynomials) we multiply them with the corresponding monomials, and hence the coefficients will become larger.  This will slow down the algorithm and increase the memory usage. In order to keep the coefficients small we have implemented the following algorithm.

\begin{alg}[preIntegerCheck]\label{alg: get constant fast}
Let $I = \gen{f_1, \ldots, f_r}$ be an ideal in $\zz[x]$ and  $<$ a fixed monomial ordering. We use the following idea to get access to  a constant or a monomial (if it exists) from $I$. This is useful when computing over $\zz$ because adding this constant to the generating system of $I$ before starting the standard basis algorithm will keep the coefficients small. 

The main idea is described as follows: we want to compute a standard basis over $\zz[x]$. We will first compute the standard basis over $\qq[x]$. If the result is 1, then we know there is a non-zero integer in the ideal.

\begin{algorithm}[H]
\caption{preIntegerCheck(ideal $I$)}
\begin{algorithmic}[1]
\REQUIRE $I = \gen{f_1, \ldots, f_r} \subseteq \zz[x]$ an ideal
\ENSURE$\aco{f_1, \ldots, f_r, c\cdot m}$, $c \in \zz$, $m$ monomial; $m = 1$ and $c \neq 0$ if $I \cap \zz \neq 0.$
\STATE $J = \aco{1,f_1, \ldots, f_r}$
\STATE compute $S$, a standard basis  of $I$ over $\qq[x]$ 
\STATE $\textnormal{compute the syzygies }Z \subset \Dsum{i=0}{r}\qq[x]\cdot \epsilon_i\textnormal{ of }J$

\IF{$S = \gen{1}$}
\STATE $\textnormal{search in }Z\textnormal{ for a syzygy where the 0 component consists of just a constant }c \cdot \epsilon_0$
\ret\ $I\cup\aco{c}$
\ELSE
\STATE $\textnormal{search in }Z\textnormal{ for a syzygy where the 0 component consists of just a term }c\cdot m \cdot \epsilon_0$
\IF{$\textnormal{such monomial is found}$}
\ret\ $I \cup \aco{c\cdot m}$
\ENDIF
\ENDIF
\ret\ $I$
\end{algorithmic}
\end{algorithm}
\end{alg}

Note that Algorithm \ref{alg: get constant fast} is very costly and if the algorithm does not find a constant or monomial, it increases the standard basis run-time without bringing any new information. However in the other cases, the algorithms run-time will improve because we have a bound on the coefficients. This is very useful over $\zz$ since the biggest problem proves to be the size of the coefficients that appear during the standard basis computation. In case of parallel systems, one could run both versions (with the preIntegerCheck and without) and stop as soon as one version finishes. 

%

\begin{ex}\label{ex 18} In this example we use the degree reverse lexicographical ordering. Let $I = \gen{x+4,xy+9,x-y+8}$ be an ideal in $\zz[x,y]$. We compute the syzygies $Z \subset  \Dsum{i=0}{3}\qq[x,y] $ of $J = \aco{1,I} =: \aco{1,f_1,f_2,f_3}$. The syzygies of $J$ are
$$\begin{array}{l}
7 \cdot \epsilon_1 - (x+4)\cdot \epsilon_2 + \epsilon_3 + x \cdot \epsilon_4, \\
(y-4) \cdot \epsilon_1 - \epsilon_2 + \epsilon_4 \textnormal{ and } \\
(x+4) \cdot \epsilon_1 - \epsilon_2.
\end{array}$$

Indeed, $$7 \cdot 1 + (-x-4) \cdot f_1 + 1 \cdot f_2 + x \cdot f_3 = 0.$$ 

Since the standard basis of $I$ over $\qq[x,y]$ is $1$, it is expected that one of the syzygies has on the first component a constant $-$ in this case it is the first syzygy. 
We continue by computing the standard basis of $\gen{7,I}$ over $\zz[x,y]$: $\aco{7,x+4,y-4}$.
\end{ex}

This was a very simple example. Now consider the ideal $I = \gen{f_1,\ldots,f_{70}} \subseteq \zz[x,y,z]$ and consider the degree reverse lexicographical ordering. The $f_i$ can be found in the Appendix.

Despite of how \textit{complicated} the ideal $I$ seems to be, its standard basis consists of few small polynomials. We give below the Gr\"obner basis $S$ of $I$:

$$
S = 
\begin{array}{l}
\big\{18,
6z-12,
2y-4,
2z^2+4z+8,
yz+z+3,
3x^2z-15x^2,
x^2y+3x^2z+x^2,\\
x^3+10z,
x^2z^2-4x^2z-11x^2
\big\}\end{array}
$$

Using the new strategy described in Algorithm \ref{alg: get constant fast}, we compute the standard basis over $\qq[x,y,z]$ and see that it is equal to $\gen{1}$. Hence the ideal contains an integer. Using $\sing$'s syzygy procedure \verb"syz", we are not lucky enough to already get access to 18, but a multiple of it: 6\,133\,248. This seems to be a big number, but if we add this integer to the generators of \verb"I" the computations will drastically speed up.

We show how 6\,133\,248 can be generated using $f_1, \ldots, f_{70}$.

$$\arraycolsep=0.1pt
\begin{array}{cr}
{{6\,133\,248}}  = & {(-\tinm{231913440}x^6z+\tinm{9838752}x^3y-\tinm{4599936}x^3z+\tinm{69696}y)} \cdot f_{1}  + \\
 & {{(-\tinm{1584}y^3+\tinm{13464}y^2+\tinm{12672}y)}}{ \cdot } { f_{2}}{ + }\\ 
 & {{(-\tinm{115956720}x^6z+\tinm{1756920}x^3y-\tinm{1724976}x^3z-\tinm{40409160}x^3+
\tinm{52272}y+\tinm{200376})}}{ \cdot } { f_{3}}{ + }\\ 
 & {{(-\tinm{38652240}x^3y^3+\tinm{38652240}x^3y^2+\tinm{154608960}x^3+\tinm{1149984}y^2
 +\tinm{4983264}y+\tinm{3066624})}}{ \cdot } { f_{4}}{ + }\\ 
 & {{(-\tinm{2551047840}x^6z+\tinm{38652240}x^3y-\tinm{37949472}x^3z-\tinm{889001520}x^3+
 \tinm{574992}y+\tinm{2779128})}}{ \cdot } { f_{5}}{ + }\\ 
 & {{(-\tinm{159720}x^3y^2-\tinm{159720}x^3y-\tinm{319440}x^3-\tinm{2376}y^2-\tinm{10296}y
 -\tinm{12672})}}{ \cdot } { f_{10}}{ + }\\ 
 & {{(-\tinm{79860}x^4+\tinm{110}xy^3+\tinm{110}xy^2+\tinm{220}xy-\tinm{2508}x)}}{ \cdot } { f_{11}}{ + }\\ 
 & {{(\tinm{26136})}}{ \cdot } { f_{13}}{ + }\\ 
 & {{(\tinm{3630}x^4y^2+\tinm{3630}x^4y+\tinm{7260}x^4+\tinm{54}xy^2+\tinm{154}xy+\tinm{96}x)}}{ \cdot } { f_{14}}{ + }\\ 
 & {{(\tinm{1540}xy+\tinm{517}x)}}{ \cdot } { f_{15}}{ + }\\ 
 & {{(-\tinm{792})}}{ \cdot } { f_{16}}{ + }\\ 
 & {{(\tinm{425920}x)}}{ \cdot } { f_{18}}{ + }\\ 
 & {{(-\tinm{115956720}x^3y^2-\tinm{115956720}x^3y-\tinm{231913440}x^3-\tinm{1724976}y^2-
 \tinm{7474896}y-\tinm{4599936})}}{ \cdot } { f_{19}}{ + }\\ 
 & {{(-\tinm{1756920}x^4+\tinm{2420}xy^3+\tinm{2420}xy^2+\tinm{4840}xy-\tinm{132616}x)}}{ \cdot } { f_{26}}{ + }\\ 
 & {{(\tinm{851840}x)}}{ \cdot } { f_{27}}{ + }\\ 
 & {{(\tinm{287496})}}{ \cdot } { f_{28}}{ + } \\
 & 
{{(-\tinm{38720}x)}}{ \cdot } { f_{38}}\phantom{+}
\end{array}$$

\section{The normal form procedure}\label{section: red local rings}

%
%
%

It is well-known that the usual normal form procedure in Buchberger's algorithm does not terminate for local orderings (cf. \cite{singular-book}, \cite{M}, \cite{MPT}). For polynomial rings over a field T. Mora solved the problems introducing a larger set T of polynomials for the reduction. For polynomial rings over $\zz$ we need not only to add the polynomials with large ecart to T but also the corresponding gcd$-$polynomials.

\begin{alg}[normal form procedure]\label{alg: reduce blueprint} The following algorithm computes the normal form over rings in case of local and mixed orderings.
\begin{algorithm}[H]
\caption{reduce($f$, $G$)}
\begin{algorithmic}[1]
\REQUIRE $f \in R[x] \textnormal{ and a finite subset }G\subset R[x], < \textnormal{ ordering}$
\ENSURE $h \in R[x], \textnormal{ a weak normal form of }f \textnormal{ w.r.t. }G, \textnormal{ i.e. } h = uf-\Sum{g \in G}{}c_gg, c_g \in R[x], u \in R[x]_< \textnormal{ unit}, \lm(f) \geq \lm(c_gg) \textnormal{ for all } g \textnormal{ such that } c_gg \neq 0 \textnormal{ and } h=0 \textnormal{ or } \lt(h) \not\in L(G).$

\STATE $h = f$
\STATE $T = G$
\WHILE{$h \neq 0$ and there is a $g \in T$ with $\lt(g)|\lt(h)$}
\STATE {choose such $g$ with minimal ecart}
\IF{$\ecart(g) > \ecart(h)$}
\STATE $T = T \cup \aco{h} \cup \set{\gp(h,t)}{t \in T}$
\ENDIF
\STATE $h = \sp(h,g)$
\ENDWHILE
\ret\ $h$
\end{algorithmic}
\end{algorithm}
\end{alg}

The proof that Algorithm \ref{alg: reduce blueprint} terminates and gives the expected result is given by \linebreak Lichtblau \cite{LB} in the case of global orderings. The termination and correctness of the algorithm for polynomial rings over a field with non-global orderings can be found in \cite{singular-book}. Combining the two proofs gives termination and correctness for Algorithm \ref{alg: reduce blueprint} in case of polynomial rings over $\zz$ and non-global orderings.

\begin{ex}\label{ex: infinite reduction}
Let $\zz[x,y]$  with the local degree reverse monomial ordering. Consider the ideal $$I = \gen{6+y+x^2, 4+x}.$$ Then a standard basis of $I$ with respect to the ordering is $$S = \aco{2-x+y+x^2,x-2y-x^2-xy-x^3}.$$ 
$$T = \aco{2-x+y+x^2,x-2y-x^2-xy-x^3},$$
\begin{equation}
 \xy
 \xymatrixcolsep{6pc}
 \xymatrix@R-2pc{
{\phantom{\boxed{^2}}} & \save {4+x}*[F]\frm{} \restore \ar@{->}`r[d]`d[d]`l[ddl]`d[dd][dd]^-{2-x+y+x^2}\\
{\phantom{\boxed{^2}}} &\\
{\phantom{\boxed{^2}}} & \save {3x-2y-2x^2}*[F]\frm{}\restore \ar@{->}`r[d]`d[d]`l[ddl]`d[dd][dd]^-{x-2y-x^2-xy+x^3} \\
{\phantom{\boxed{^2}}} &\\ 
{\phantom{\boxed{^2}}} & 4y+x^2+3xy+3x^3 \ar@{->}`r[d]`d[d]`l[ddl]`d[dd][dd]^-{2-x+y+x^2}\\
{\phantom{\boxed{^2}}} &\\
{\phantom{\boxed{^2}}} & x^2+5xy-2y^2+3x^3-2x^2y \ar@{-->}`r[d]`d[d]`l[ddl]`d[dd][dd]^-{x+2y+3x^2}\\
{\phantom{\boxed{^2}}} & \\
{\phantom{\boxed{^2}}}  & 3xy-2y^2-2x^2y\ar@{-->}`r[d]`d[d]`l[ddl]`d[dd][dd]^-{3x-2y-2x^2} \\
{\phantom{\boxed{^2}}} &\\
{\phantom{\boxed{^2}}} & 0 \\
}
\endxy
\xy\xymatrix@R-2pc{{\phantom{\boxed{^2}}}\\ {\phantom{\boxed{^2}}}\\
T = T \cup \aco{\boxed{x+2y+3x^2}}
}
\endxy
\end{equation}
At the end of the reduction $$T = \aco{2-x+y+x^2,x-2y-x^2-xy-x^3,4+x,3x-2y+3x^2,x+2y+3x^2}.$$

Note that Algorithm \ref{alg: reduce blueprint} does not terminate in the example above if the gcd$-$polynoamial $x+2y+3x^2$ would not be added to T.
\end{ex}



\section{Appendix}

In the following tables we print the ideals corresponding to the examples presented when comparing the ALL vs JUST timings in Section \ref{section: all vs just}.
\label{appendix: all vs just}
\begin{table}[H]
\begin{center}
\arraycolsep=2.4pt\def\arraystretch{1}
\begin{tabular}{l||l||}
 & \multicolumn{1}{c|}{Generators for ideal}  \\ \hline \hline
Example \randex{A}{1} \hide{2:106} &  $\begin{array}{l}\tinm{17} \cdot x^2y,\\

\tinm{22} \cdot y^3z+\tinm{3} \cdot x^2z^2+\tinm{28} \cdot y^2z^2+\tinm{9} \cdot yz^2+\tinm{83} \cdot x^2+\tinm{13} \cdot yz,\\

\tinm{66} \cdot y^3z+\tinm{63} \cdot xyz^2+\tinm{85} \cdot z^3\end{array}$
 \\ \hline

Example \randex{A}{2} \hide{1:250} & $\begin{array}{l}   \tinm{27} \cdot xyz+\tinm{13} \cdot x^2+\tinm{89} \cdot y^2+\tinm{42} \cdot xz,\\
\tinm{35} \cdot x^3+\tinm{68} \cdot xy,\\
\tinm{44} \cdot x^3+\tinm{13} \cdot y^3+\tinm{81} \cdot y^2z+\tinm{4} \cdot yz^2\\  \end{array}$  \\ \hline

Example \randex{A}{3} \hide{2:109} & $\begin{array}{l}  \tinm{98} \cdot x^3y+\tinm{45} \cdot yz^3,\\
\tinm{16} \cdot xy^3+\tinm{50} \cdot x^2y+\tinm{45} \cdot y^2+\tinm{82} \cdot z^2,\\
\tinm{9} \cdot x^3y+\tinm{49} \cdot x^2yz+\tinm{61} \cdot y^2z+\tinm{52} \cdot z^3   \end{array}$  \\ \hline

Example \randex{A}{4} \hide{1:238} & $\begin{array}{l}  \tinm{25} \cdot y^3,\\
\tinm{60} \cdot y^3+\tinm{12} \cdot xyz+\tinm{54} \cdot y^2z+\tinm{98} \cdot yz^2+\tinm{35} \cdot x^2+\tinm{88} \cdot xy+\tinm{19} \cdot z^2,\\
\tinm{87} \cdot x^2y+\tinm{96} \cdot x^2   \end{array}$ \\ \hline

Example \randex{A}{5} \hide{1:227} & $\begin{array}{l}  \tinm{76} \cdot y^2z+\tinm{61} \cdot y^2+\tinm{51} \cdot yz+\tinm{19} \cdot z^2,\\
\tinm{31} \cdot x^3+\tinm{3} \cdot xy^2+\tinm{70} \cdot y^2z,\\
\tinm{84} \cdot x^3+\tinm{30} \cdot x^2z+\tinm{44} \cdot xz   \end{array}$ \\ \hline 

Example \randex{A}{6} \hide{1:265} & $\begin{array}{l}  \tinm{19} \cdot x^3+\tinm{2} \cdot xy+\tinm{29} \cdot y^2,\\
\tinm{9} \cdot xy^2+\tinm{42} \cdot y^3+\tinm{2} \cdot yz^2,\\
\tinm{54} \cdot xy^2+\tinm{83} \cdot x^2z+\tinm{98} \cdot xy+\tinm{78} \cdot yz   \end{array}$ \\ \hline 

Example \randex{A}{7} \hide{2:125} & $\begin{array}{l}  \tinm{72} \cdot x^3y+\tinm{50} \cdot x^2yz,\\
\tinm{64} \cdot x^3z+\tinm{30} \cdot x^2yz+\tinm{74} \cdot x^3+\tinm{38} \cdot xy^2+\tinm{74} \cdot z^3,\\
\tinm{76} \cdot x^2y^2+\tinm{13} \cdot y^2+\tinm{40} \cdot z^2   \end{array}$ \\ \hline

Example \randex{A}{8} \hide{2:105} & $\begin{array}{l}  \tinm{27} \cdot xz^3+\tinm{87} \cdot y^2z+\tinm{3} \cdot z^2,\\
\tinm{67} \cdot x^3z+\tinm{42} \cdot y^3z+\tinm{67} \cdot x^2y+\tinm{90} \cdot yz^2+\tinm{73} \cdot xy,\\
\tinm{38} \cdot x^3y+\tinm{69} \cdot x^2yz   \end{array}$ \\ \hline 

Example \randex{A}{9} \hide{2:141} & $\begin{array}{l}  \tinm{21} \cdot x^3y+\tinm{18} \cdot x^3z+\tinm{45} \cdot x^2yz+\tinm{100} \cdot xy^2z+\tinm{43} \cdot yz,\\
\tinm{85} \cdot x^2y+\tinm{93} \cdot xy^2,\\
\tinm{14} \cdot y^2z^2+\tinm{6} \cdot x^2y+\tinm{91} \cdot z^2   \end{array}$ \\ \hline 

Example \randex{A}{10} \hide{0:66} & $\begin{array}{l}  \tinm{69} \cdot y^2,\\
\tinm{47} \cdot x^3y+\tinm{82} \cdot xyz^2+\tinm{74} \cdot yz^3+\tinm{55} \cdot xyz+\tinm{96} \cdot xz^2+\tinm{46} \cdot x^2,\\
\tinm{16} \cdot xy^2z+\tinm{17} \cdot z^4+\tinm{36} \cdot yz^2   \end{array}$ \\ \hline \hline 
\end{tabular}
\end{center}
\caption{Examples for ALL Strategy}
\end{table}
\begin{table}[H]
\begin{center}
\arraycolsep=2.4pt\def\arraystretch{1.3}
\begin{tabular}{l||l||}
 & \multicolumn{1}{c|}{Generators for ideal}  \\ \hline \hline

Example \randex{B}{1} \hide{3:1} & $\begin{array}{l}  \tinm{6} \cdot x^3z+\tinm{29} \cdot x^2z+\tinm{42} \cdot xy,\\
x^3z+\tinm{47} \cdot x^2yz+\tinm{28} \cdot xz^2+\tinm{46} \cdot x^2,\\
\tinm{96} \cdot z^3   \end{array}$ \\ \hline

Example \randex{B}{2} \hide{0:79} & $\begin{array}{l}  \tinm{9} \cdot x^2y^2+\tinm{51} \cdot x^3z+\tinm{10} \cdot z^3+\tinm{28} \cdot x^2+\tinm{7} \cdot y^2,\\
\tinm{43} \cdot x^3y+\tinm{3} \cdot x^2z^2+\tinm{86} \cdot xyz^2+\tinm{24} \cdot z^4+\tinm{67} \cdot x^2z+\tinm{68} \cdot yz^2+\tinm{27} \cdot xy,\\
\tinm{23} \cdot xz^2   \end{array}$  \\ \hline

Example \randex{B}{3} \hide{3:8} & $\begin{array}{l}  \tinm{50} \cdot xz^3+\tinm{49} \cdot yz^2+\tinm{15} \cdot z^2,\\
\tinm{2} \cdot x^3y+\tinm{16} \cdot y^3z+\tinm{74} \cdot y^3+\tinm{53} \cdot x^2,\\
\tinm{4} \cdot xy^2   \end{array}$ \\ \hline

Example \randex{B}{4} \hide{3:6} & $\begin{array}{l}  \tinm{57} \cdot xyz^2+\tinm{32} \cdot y^3+\tinm{26} \cdot yz^2+\tinm{24} \cdot z^2,\\
\tinm{27} \cdot y^4+\tinm{33} \cdot y^3z+\tinm{94} \cdot z^2,\\
\tinm{52} \cdot x^2z   \end{array}$ \\ \hline

Example \randex{B}{5} \hide{3:12} & $\begin{array}{l}   \tinm{91} \cdot y^4+\tinm{20} \cdot x^3z+\tinm{34} \cdot x^2,\\
\tinm{38} \cdot xyz^2+\tinm{18} \cdot x^3+\tinm{95} \cdot x^2z+\tinm{82} \cdot yz,\\
\tinm{98} \cdot y^3  \end{array}$ \\ \hline \hline

\end{tabular}
\end{center}
\caption{Examples for JUST Strategy}
\end{table}


\label{appendix: example 18}
The ideal for Example \ref{ex 18}: 

We consider over $\zz[x,y,z]$ the degree reverse lexicographical ordering and the ideal \verb"I" generated by the following 70 polynomials $f_1,\ldots, f_{70}$.

$$\begin{array}{ccl}
f_{1} & = &\tinm{42} \cdot x^3z+y^2z-yz+\tinm{11} \cdot y\\
f_{2} & = &y^3z^2-y^2z^2+\tinm{11} \cdot y^2z+\tinm{484}\\
f_{3} & = &y^4z-y^3z-\tinm{10648} \cdot x^3+\tinm{11} \cdot y^3-\tinm{44} \cdot y^2+\tinm{44} \cdot y\\
f_{4} & = &x^3yz^2-\tinm{2}\\
f_{5} & = &\tinm{11} \cdot x^3y^2z+\tinm{484} \cdot x^3+\tinm{2} \cdot y^2-\tinm{2} \cdot y\\
f_{6} & = &\tinm{117128} \cdot x^6-\tinm{121} \cdot x^3y^3+\tinm{968} \cdot x^3y^2-\tinm{968} \cdot x^3y+\tinm{2} \cdot y^4-\tinm{4} \cdot y^3+\tinm{2} \cdot y^2\\
f_{7} & = &\tinm{121} \cdot x^6z^3+yz-z+\tinm{11}\\
f_{8} & = &\tinm{2178} \cdot x^2y^2z^3-\tinm{1452} \cdot x^2yz^3+\tinm{15972} \cdot x^2yz^2\\
f_{9} & = &\tinm{1452} \cdot x^2y^3z^2-\tinm{1452} \cdot x^2y^2z^2+\tinm{7986} \cdot x^2y^2z\\
f_{10} & = &-\tinm{726} \cdot x^3y^2z^2+\tinm{484} \cdot x^3yz^2+y^4z^2-\tinm{5324} \cdot x^3yz-y^3z^2+\tinm{22} \cdot y^3z-\tinm{11} \cdot y^2z\\&&+\tinm{121} \cdot y^2\\

f_{11} & = &\tinm{2904} \cdot x^2y^3z^2-\tinm{2178} \cdot x^2y^2z^2+\tinm{23958} \cdot x^2y^2z+\tinm{351384} \cdot x^2\\
f_{12} & = &\tinm{726} \cdot x^2y^4z-\tinm{726} \cdot x^2y^3z+\tinm{7730448} \cdot x^5+\tinm{31944} \cdot x^2y^2-\tinm{31944} \cdot x^2y\\
\end{array}$$
$$\begin{array}{ccl}
f_{13} & = &-\tinm{968} \cdot x^3y^3z+\tinm{726} \cdot x^3y^2z-\tinm{2} \cdot y^5z-\tinm{7986} \cdot x^3y^2+\tinm{4} \cdot y^4z+\tinm{21296} \cdot x^3y\\&&-\tinm{22} \cdot y^4-\tinm{2} \cdot y^3z-\tinm{10648} \cdot x^3+\tinm{110} \cdot y^3-\tinm{132} \cdot y^2+\tinm{44} \cdot y\\

f_{14} & = &\tinm{95832} \cdot x^2y^2z^2-\tinm{63888} \cdot x^2yz^2+\tinm{702768} \cdot x^2yz\\
f_{15} & = &\tinm{63888} \cdot x^2y^3z-\tinm{63888} \cdot x^2y^2z+\tinm{351384} \cdot x^2y^2\\

f_{16} & = &-\tinm{5} \cdot y^6z^2+\tinm{9} \cdot y^5z^2-\tinm{88} \cdot y^5z-\tinm{4} \cdot y^4z^2+\tinm{253} \cdot y^4z-\tinm{363} \cdot y^4-\tinm{264} \cdot y^3z\\&&+\tinm{968} \cdot y^3+\tinm{88} \cdot y^2z-\tinm{484} \cdot y^2\\
f_{17} & = &\tinm{726} \cdot x^5z^3-\tinm{6} \cdot x^2y^2z^3+\tinm{3} \cdot x^2yz^3-\tinm{33} \cdot x^2yz^2\\
f_{18} & = &\tinm{726} \cdot x^5yz^2-\tinm{3} \cdot x^2y^3z^2+\tinm{3} \cdot x^2y^2z^2\\
f_{19} & = &-\tinm{242} \cdot x^6z^2+\tinm{3} \cdot x^3y^2z^2-x^3yz^2+\tinm{22} \cdot x^3yz\\
f_{20} & = &-\tinm{9} \cdot x^2y^3z^4+\tinm{6} \cdot x^2y^2z^4-\tinm{66} \cdot x^2y^2z^3\\
f_{21} & = &-\tinm{6} \cdot x^2y^4z^3+\tinm{6} \cdot x^2y^3z^3-\tinm{33} \cdot x^2y^3z^2\\
f_{22} & = &\tinm{4} \cdot x^3y^3z^3-\tinm{2} \cdot x^3y^2z^3+\tinm{33} \cdot x^3y^2z^2\\
f_{23} & = &-\tinm{12} \cdot x^2y^4z^3+\tinm{9} \cdot x^2y^3z^3-\tinm{31944} \cdot x^5z^2-\tinm{99} \cdot x^2y^3z^2+\tinm{264} \cdot x^2y^2z^2-\tinm{132} \cdot x^2yz^2\\
f_{24} & = &-\tinm{3} \cdot x^2y^5z^2+\tinm{3} \cdot x^2y^4z^2-\tinm{63888} \cdot x^5yz\\
f_{25} & = &\tinm{7} \cdot x^3y^4z^2-\tinm{5} \cdot x^3y^3z^2+\tinm{66} \cdot x^3y^3z-\tinm{176} \cdot x^3y^2z+\tinm{88} \cdot x^3yz\\
f_{26} & = &\tinm{15972} \cdot x^5yz^2-\tinm{66} \cdot x^2y^3z^2+\tinm{33} \cdot x^2y^2z^2-\tinm{363} \cdot x^2y^2z-\tinm{15972} \cdot x^2\\
f_{27} & = &-\tinm{33} \cdot x^2y^4z+\tinm{33} \cdot x^2y^3z-\tinm{351384} \cdot x^5-\tinm{1452} \cdot x^2y^2+\tinm{1452} \cdot x^2y\\

f_{28} & = &-\tinm{5324} \cdot x^6yz+\tinm{11} \cdot x^3y^2z+\tinm{121} \cdot x^3y^2-\tinm{968} \cdot x^3y+\tinm{484} \cdot x^3-\tinm{4} \cdot y^3+\tinm{6} \cdot y^2-\tinm{2} \cdot y\\
f_{29} & = &-\tinm{99} \cdot x^2y^4z^3+\tinm{66} \cdot x^2y^3z^3-\tinm{726} \cdot x^2y^3z^2-\tinm{4356} \cdot x^2y^2z^2+\tinm{2904} \cdot x^2yz^2\\&&-\tinm{31944} \cdot x^2yz\\
f_{30} & = &-\tinm{66} \cdot x^2y^5z^2+\tinm{66} \cdot x^2y^4z^2-\tinm{363} \cdot x^2y^4z-\tinm{2904} \cdot x^2y^3z+\tinm{2904} \cdot x^2y^2z\\&&-\tinm{15972} \cdot x^2y^2\\
f_{31} & = &-\tinm{11} \cdot x^3y^4z^2+\tinm{22} \cdot x^3y^3z^2-\tinm{8} \cdot y^4z+\tinm{12} \cdot y^3z-\tinm{44} \cdot y^3-\tinm{4} \cdot y^2z+\tinm{22} \cdot y^2\\
f_{32} & = &-\tinm{132} \cdot x^2y^5z^2+\tinm{99} \cdot x^2y^4z^2-\tinm{702768} \cdot x^5yz-\tinm{1089} \cdot x^2y^4z-\tinm{2904} \cdot x^2y^3z\\&&+\tinm{2904} \cdot x^2y^2z-\tinm{47916} \cdot x^2y^2\\
f_{33} & = &-\tinm{33} \cdot x^2y^6z+\tinm{33} \cdot x^2y^5z-\tinm{351384} \cdot x^5y^2-\tinm{1452} \cdot x^2y^4+\tinm{1452} \cdot x^2y^3\\
f_{34} & = &\tinm{22} \cdot x^3y^5z-\tinm{11} \cdot x^3y^4z+\tinm{363} \cdot x^3y^4-\tinm{968} \cdot x^3y^3+\tinm{484} \cdot x^3y^2-\tinm{4} \cdot y^5+\tinm{6} \cdot y^4\\&&-\tinm{2} \cdot y^3\\
f_{35} & = &\tinm{33} \cdot x^5y^2z^3-\tinm{1452} \cdot x^5z^2+\tinm{12} \cdot x^2y^2z^2-\tinm{6} \cdot x^2yz^2\\
f_{36} & = &-\tinm{33} \cdot x^5y^3z^2-\tinm{2904} \cdot x^5yz\\
f_{37} & = &-\tinm{33} \cdot x^6y^2z^2-\tinm{8} \cdot x^3y^2z+\tinm{4} \cdot x^3yz\\
f_{38} & = &-\tinm{263538} \cdot x^5y^2z+\tinm{726} \cdot x^2y^4z-\tinm{363} \cdot x^2y^3z-\tinm{7730448} \cdot x^5+\tinm{3993} \cdot x^2y^3\\&&-\tinm{31944} \cdot x^2y^2+\tinm{31944} \cdot x^2y\\
\end{array}$$
$$\begin{array}{ccl}
f_{39} & = &-\tinm{170069856} \cdot x^8+\tinm{87846} \cdot x^5y^3-\tinm{1405536} \cdot x^5y^2+\tinm{363} \cdot x^2y^5+\tinm{1405536} \cdot x^5y\\&&-\tinm{3267} \cdot x^2y^4+\tinm{5808} \cdot x^2y^3-\tinm{2904} \cdot x^2y^2\\
f_{40} & = &\tinm{87846} \cdot x^6y^2-\tinm{468512} \cdot x^6y+\tinm{363} \cdot x^3y^4+\tinm{234256} \cdot x^6-\tinm{4235} \cdot x^3y^3+\tinm{5808} \cdot x^3y^2\\&&-\tinm{8} \cdot y^5-\tinm{1936} \cdot x^3y+\tinm{20} \cdot y^4-\tinm{16} \cdot y^3+\tinm{4} \cdot y^2\\
f_{41} & = &-\tinm{2108304} \cdot x^5y^2z^2+\tinm{1089} \cdot x^2y^5z^2+\tinm{1405536} \cdot x^5yz^2-\tinm{9438} \cdot x^2y^4z^2\\&&-\tinm{15460896} \cdot x^5yz+\tinm{7986} \cdot x^2y^4z+\tinm{14520} \cdot x^2y^3z^2\\&&-\tinm{63888} \cdot x^2y^3z-\tinm{5808} \cdot x^2y^2z^2+\tinm{63888} \cdot x^2y^2z\\
f_{42} & = &-\tinm{1405536} \cdot x^5y^3z+\tinm{726} \cdot x^2y^6z+\tinm{1405536} \cdot x^5y^2z-\tinm{6534} \cdot x^2y^5z-\tinm{7730448} \cdot x^5y^2\\&&+\tinm{3993} \cdot x^2y^5+\tinm{11616} \cdot x^2y^4z-\tinm{31944} \cdot x^2y^4-\tinm{5808} \cdot x^2y^3z+\tinm{31944} \cdot x^2y^3\\
f_{43} & = &\tinm{726} \cdot x^3y^5z-\tinm{4598} \cdot x^3y^4z+\tinm{3993} \cdot x^3y^4+\tinm{5808} \cdot x^3y^3z-\tinm{16} \cdot y^6z-\tinm{21296} \cdot x^3y^3\\&&-\tinm{1936} \cdot x^3y^2z+\tinm{40} \cdot y^5z+\tinm{10648} \cdot x^3y^2-\tinm{88} \cdot y^5-\tinm{32} \cdot y^4z+\tinm{132} \cdot y^4+\tinm{8} \cdot y^3z\\&&-\tinm{44} \cdot y^3\\

f_{44} & = &-\tinm{2811072} \cdot x^5y^3z+\tinm{1452} \cdot x^2y^6z+\tinm{2108304} \cdot x^5y^2z-\tinm{12705} \cdot x^2y^5z-\tinm{11595672} \cdot x^5y^2\\&&+\tinm{11979} \cdot x^2y^5+\tinm{20328} \cdot x^2y^4z-\tinm{127776} \cdot x^2y^4-\tinm{8712} \cdot x^2y^3z+\tinm{111804} \cdot x^2y^3\\
f_{45} & = &-\tinm{702768} \cdot x^5y^4+\tinm{363} \cdot x^2y^7+\tinm{702768} \cdot x^5y^3-\tinm{3267} \cdot x^2y^6+\tinm{5808} \cdot x^2y^5-\tinm{2904} \cdot x^2y^4\\
f_{46} & = &\tinm{363} \cdot x^3y^6-\tinm{2299} \cdot x^3y^5+\tinm{2904} \cdot x^3y^4-\tinm{8} \cdot y^7-\tinm{968} \cdot x^3y^3+\tinm{20} \cdot y^6-\tinm{16} \cdot y^5+\tinm{4} \cdot y^4\\
f_{47} & = &-\tinm{702768} \cdot x^8z^2-\tinm{726} \cdot x^5y^3z^2+\tinm{2904} \cdot x^5y^2z^2+\tinm{24} \cdot x^2y^4z^2-\tinm{36} \cdot x^2y^3z^2\\&&+\tinm{12} \cdot x^2y^2z^2\\
f_{48} & = &-\tinm{1405536} \cdot x^8yz+\tinm{726} \cdot x^5y^4z-\tinm{5808} \cdot x^5y^3z+\tinm{5808} \cdot x^5y^2z\\
f_{49} & = &\tinm{726} \cdot x^6y^3z-\tinm{3872} \cdot x^6y^2z+\tinm{1936} \cdot x^6yz-\tinm{16} \cdot x^3y^4z+\tinm{24} \cdot x^3y^3z-\tinm{8} \cdot x^3y^2z\\
f_{50} & = &-\tinm{15460896} \cdot x^8yz-\tinm{3993} \cdot x^5y^4z+\tinm{31944} \cdot x^5y^2z+\tinm{264} \cdot x^2y^5z-\tinm{527076} \cdot x^5y^2\\&&-\tinm{396} \cdot x^2y^4z+\tinm{1452} \cdot x^2y^4+\tinm{132} \cdot x^2y^3z-\tinm{726} \cdot x^2y^3\\
f_{51} & = &-\tinm{7730448} \cdot x^8y^2+\tinm{3993} \cdot x^5y^5-\tinm{31944} \cdot x^5y^4+\tinm{31944} \cdot x^5y^3\\
f_{52} & = &\tinm{3993} \cdot x^6y^4-\tinm{21296} \cdot x^6y^3+\tinm{10648} \cdot x^6y^2-\tinm{88} \cdot x^3y^5+\tinm{132} \cdot x^3y^4-\tinm{44} \cdot x^3y^3\\
f_{53} & = &-\tinm{1452} \cdot x^5yz^4+\tinm{726} \cdot x^5z^4-\tinm{7986} \cdot x^5z^3+\tinm{726} \cdot x^2z^2\\
f_{54} & = &\tinm{87846} \cdot x^8z^3-\tinm{726} \cdot x^5y^2z^3+\tinm{726} \cdot x^5yz^3+\tinm{726} \cdot x^2yz-\tinm{726} \cdot x^2z\\
f_{55} & = &\tinm{726} \cdot x^6yz^3-\tinm{363} \cdot x^6z^3+\tinm{3993} \cdot x^6z^2-\tinm{242} \cdot x^3z+y^2z-\tinm{2} \cdot yz+\tinm{11} \cdot y+z-\tinm{11}\\
f_{56} & = &-\tinm{2178} \cdot x^5y^2z^5+\tinm{1452} \cdot x^5yz^5-\tinm{15972} \cdot x^5yz^4\\

f_{57} & = &-\tinm{1452} \cdot x^5y^3z^4+\tinm{1452} \cdot x^5y^2z^4-\tinm{7986} \cdot x^5y^2z^3\\
f_{58} & = &\tinm{1089} \cdot x^6y^2z^4-\tinm{726} \cdot x^6yz^4+\tinm{7986} \cdot x^6yz^3+y^3z^2-\tinm{3} \cdot y^2z^2+\tinm{11} \cdot y^2z+\tinm{2} \cdot yz^2\\&&-\tinm{22} \cdot yz\\
f_{59} & = &-\tinm{2904} \cdot x^5y^3z^4+\tinm{2178} \cdot x^5y^2z^4-\tinm{23958} \cdot x^5y^2z^3+\tinm{63888} \cdot x^5yz^3-\tinm{31944} \cdot x^5z^3\\&&-\tinm{31944} \cdot x^2z\\
\end{array}$$
$$\begin{array}{ccl}
f_{60} & = &-\tinm{726} \cdot x^5y^4z^3+\tinm{726} \cdot x^5y^3z^3-\tinm{11595672} \cdot x^8z^2-\tinm{31944} \cdot x^2y+\tinm{31944} \cdot x^2\\
f_{61} & = &\tinm{1452} \cdot x^6y^3z^3-\tinm{1089} \cdot x^6y^2z^3+\tinm{11979} \cdot x^6y^2z^2-\tinm{31944} \cdot x^6yz^2+\tinm{15972} \cdot x^6z^2\\&&+\tinm{3} \cdot y^4z-\tinm{6} \cdot y^3z+\tinm{33} \cdot y^3+\tinm{3} \cdot y^2z-\tinm{121} \cdot y^2+\tinm{132} \cdot y-\tinm{44}\\

f_{62} & = &-\tinm{726} \cdot x^8z^5+\tinm{3} \cdot x^2yz^3\\
f_{63} & = &-\tinm{363} \cdot x^8yz^4+\tinm{3} \cdot x^2y^2z^2-\tinm{3} \cdot x^2yz^2\\
f_{64} & = &\tinm{363} \cdot x^9z^4-x^3yz^2-x^3z^2\\
f_{65} & = &-\tinm{15972} \cdot x^8yz^4-\tinm{2904} \cdot x^5yz^3+\tinm{1452} \cdot x^5z^3+\tinm{33} \cdot x^2y^2z^2+\tinm{1452} \cdot x^2z\\
f_{66} & = &\tinm{3993} \cdot x^8y^2z^3+\tinm{527076} \cdot x^8z^2+\tinm{33} \cdot x^2y^3z-\tinm{33} \cdot x^2y^2z+\tinm{1452} \cdot x^2y-\tinm{1452} \cdot x^2\\
f_{67} & = &\tinm{7986} \cdot x^9yz^3+\tinm{1452} \cdot x^6yz^2-\tinm{726} \cdot x^6z^2+\tinm{11} \cdot x^3y^2z-\tinm{22} \cdot x^3yz+\tinm{4} \cdot y^2-\tinm{6} \cdot y+\tinm{2}\\
f_{68} & = &\tinm{263538} \cdot x^8y^2z^3-\tinm{1405536} \cdot x^8yz^3+\tinm{702768} \cdot x^8z^3-\tinm{5808} \cdot x^5y^3z^3+\tinm{8712} \cdot x^5y^2z^3\\&&-\tinm{2904} \cdot x^5yz^3+\tinm{702768} \cdot x^5z-\tinm{363} \cdot x^2y^3z+\tinm{2904} \cdot x^2y^2z-\tinm{2904} \cdot x^2yz\\
f_{69} & = &\tinm{255104784} \cdot x^11z^2-\tinm{131769} \cdot x^8y^3z^2+\tinm{1054152} \cdot x^8y^2z^2-\tinm{1054152} \cdot x^8yz^2\\&&+\tinm{702768} \cdot x^5y-\tinm{363} \cdot x^2y^4-\tinm{702768} \cdot x^5+\tinm{3267} \cdot x^2y^3-\tinm{5808} \cdot x^2y^2+\tinm{2904} \cdot x^2y\\
f_{70} & = &-\tinm{131769} \cdot x^9y^2z^2+\tinm{702768} \cdot x^9yz^2-\tinm{351384} \cdot x^9z^2+\tinm{2904} \cdot x^6y^3z^2-\tinm{4356} \cdot x^6y^2z^2\\&&+\tinm{1452} \cdot x^6yz^2-\tinm{363} \cdot x^3y^3+\tinm{2299} \cdot x^3y^2-\tinm{2904} \cdot x^3y+\tinm{8} \cdot y^4+\tinm{968} \cdot x^3-\tinm{20} \cdot y^3\\&&+\tinm{16} \cdot y^2-\tinm{4} \cdot y\\
\end{array}$$



\end{document}